\documentclass[18pt]{article}
\usepackage{latexsym,amsfonts,amssymb,amsmath}
\usepackage[cp1251]{inputenc}
\usepackage[ukrainian, russian, english]{babel}
\usepackage[dvips]{graphicx}
\usepackage[usenames]{color}

\numberwithin{equation}{section}

\newtheorem{theorem}{Theorem}[section]
\newtheorem{lemma}{Lemma}[section]

\newtheorem{definition}{Definition}[section]

\numberwithin{equation}{section}

\newcommand{\conv}{\operatorname{conv\,}}
\newcommand{\diam}{\operatorname{diam\,}}

\newcommand{\dl}{\delta}
\newcommand{\Dl}{\Delta}
\newcommand{\e}{\varepsilon}

\newcommand{\lm}{\lambda}

\newcommand{\s}{\sigma}

\newcommand{\intt}{\operatorname{int\,}}

\newcommand{\wt}{\widetilde}
\newcommand{\Or}{\operatorname{Or\,}}

\begin{document}
\selectlanguage{english} \thispagestyle{empty}

 \pagestyle{myheadings}

\begin{center}

ON  THE STABILITY OF SOLUTIONS OF CERTAIN LINEAR  SET DIFFERENTIAL EQUATIONS

V.I. Slyn'ko

S.P. Timoshenko Institute of Mechanics of NAS of Ukraine, Kiev, Ukraine

E-mail: vitstab@ukr.net
\end{center}

\begin{center}
\textbf{Abstract}
\end{center}

New approaches to the study of  stability of solutions  of Set Differential Equations (SDEs)  based on convex geometry and the theory of mixed volumes are proposed. The stability of the forms of program solutions of linear SDEs with a stable operator is proved.  We consider the orbit of  the  action of homotheties group on the space   of  nonempty convex compacts ($\conv\Bbb R^n$)  as the form of a convex compact. For equations with periodic operator in the two-dimensional space  the asymptotic stability conditions are established.

\bigskip
\textbf{Keywords:}  Set Differential Equations, comparison method, direct Lyapunov method, Brunn–Minkowski inequality, Lyapunov stability, stable operator
\bigskip

\textbf{ MSC:} 93D30; 52A39
\section{Introduction}
Differential equations with  Hukuhara derivative (Set Differential Equations (SDEs))  were first considered  in \cite{blasi}. Further development of the theory of differential equations with  Hukuhara derivative  has been summarized in the monograph \cite{lak}, where  the conditions of existence and uniqueness of solutions of the Cauchy problem, the convergence of successive approximations including  the principle of comparison and the theorems of Lyapunov's direct method have been formulated. In papers \cite{sl-1,sl-2}  the results and  methods  of geometry  of
convex bodies, developed in the classical works of H. Minkowski and A.D. Alexandrov \cite{alex, bon}, were used for the study of stability  of solutions for dynamical systems  in the space of convex compact sets in $\Bbb R^n$. In this paper, these ideas are applied to the study of stability of solutions for Set Differential Equations.

Let $(\conv\Bbb R^n,d_H)$ be a metric space of  nonempty convex compact sets in  $\Bbb R^n$ and  $d_H$  is the Hausdorff metric.

 This work is devoted to the study of the properties of solutions for SDEs   of the form
\begin{equation}\label{21}
\gathered
D_HX(t)=\bold AX,
\endgathered
\end{equation}
where  $X(t)\in\conv\Bbb R^n$, $\bold A\in L(\Bbb R^n)$. Here and later, if $(X,\|.\|_X)$ is a Banach space, then
$L(X)$ is a Banach algebra of bounded linear operators  on    $X$.

We note that the dynamic properties of the differential equation \eqref {21}  are significantly different from the properties  of the similar ordinary differential equation (ODE)
\begin{equation}\label{22}
\gathered
\frac{dx}{dt}=\bold Ax,
\endgathered
\end{equation}
where $x\in\Bbb R^n$, $\bold A\in L(\Bbb R^n)$.

Consider the following simple example \cite{lak}. Consider the ODE
\begin{equation}\label{23}
\gathered
\frac{dx}{dt}=-x,
\endgathered
\end{equation}
$x\in\Bbb R$. A similar equation in  $\conv\Bbb R$ is of the form
\begin{equation}\label{24}
\gathered
D_HX(t)=\mathcal JX(t),
\endgathered
\end{equation}
where $X\in\conv\Bbb R$, $\mathcal J$ is a reflection operator, i.e.,
\begin{equation*}
\gathered
\mathcal JX=\{-x\,\,|\,\,x\in X\}.
\endgathered
\end{equation*}
Let $X(t)=[x_1(t),x_2(t)]\in\conv\Bbb R$, $x_1(t)\le x_2(t)$ is a solution of differential equation \eqref{24} with the initial condition
$X_0=[x_{10},x_{20}]$, $x_{10}\le x_{20}$. Then the Cauchy problem for  \eqref{24} is equivalent to  Cauchy problem for the system of differential equations
\begin{equation*}
\gathered
\frac{dx_1}{dt}=-x_2,\quad x_1(0)=x_{10},\quad \frac{dx_2}{dt}=-x_1,\quad x_2(0)=x_{20}.
\endgathered
\end{equation*}

By integrating this system we obtain the solution of the Cauchy problem \eqref{24}
\begin{equation*}
\gathered
X(t)=[x_{10}\cosh t-x_{20}\sinh t,x_{20}\cosh t-x_{10}\sinh t],\quad t\ge 0.
\endgathered
\end{equation*}
We note that $\diam X(t)=e^t\diam X_0$ and  the solution $x=0$ of ODE \eqref{23} is asymptotically stable in the sense of Lyapunov.

Consider the  stability  problem of the solution $ X \equiv 0 $  of  \eqref{24} with respect to the Hausdorff metric. It is easy to see that $d_H(X,\theta)=\max[|x_1|,|x_2|]$ and
$d_H(X(t),\theta)\to\infty$ for $t\to\infty$, provided that $ X (0) = X_0 $ is not  a single point. Thus, the solution $ X = 0 $ is unstable.

The use of Hausdorff metric as a measure does not lead to a meaningful  problem statement  about stability of solutions of SDEs. This is due to the fact that the nondecreasing  of function $\diam X(t)$ is  necessary condition for  Hukuhara differentiability of mapping $ X (t) $.

This example shows that a meaningful  problem statement about stability of solutions of SDEs
is a non-trivial task and is concerned with an adequate choice of measures with respect to which we  can consider the stability problem.

In this paper, the problem of choosing an appropriate measures with respect to which  we investigate stability is  solved on the basis of geometrical considerations.

We introduce the space of shapes of convex bodies as the quotient set of the space $\conv\Bbb R^n$  on   action  of   group of homotheties  in space $\Bbb R^n$.
Then in a conventional manner  the quotient metric for  Hausdorff metric $ d_H $ is introduced. The geometric meaning of quotient metric is the deviation of  the convex bodies shapes.  The  stability problem  is considered with respect to this quotient metric.

We note that stability of shapes of attainable sets for linear impulsive systems was considered in \cite{ovseev}. Here, the space of shape of convex compact sets  is considered as a quotient space of $ \conv \Bbb R^n $ on the action of the general affine group.

 \section{Problem statement}  Let $\mathfrak{G}$ be a certain affine group in the space $ \Bbb R^ n $, then its action  naturally extends to the space  $\conv\Bbb R^n$
\begin{equation*}
\mathfrak{g}X=\{\mathfrak{g}x\,|\,x\in X\}, \quad X\in\conv\Bbb R^n, \mathfrak{g}\in\mathfrak{G}.
\end{equation*}
The orbit of $X$, under the action of the group $\mathfrak{G}$ is defined as a subset
\begin{equation*}
\Or_{\mathfrak{G}}(X)=\{\mathfrak{g}X\,|\,\mathfrak{g}\in\mathfrak{G}\}\subset \conv\Bbb R^n.
\end{equation*}
The set of all orbits is denoted by $\conv\Bbb R^n/\mathfrak{G}$. This quotient space is endorsed with the following metric
\begin{equation*}
\rho[\Or_{\mathfrak{G}}(X),\Or_{\mathfrak{G}}(Y)]=\inf\{d_H(\mathfrak{g_1}X,\mathfrak{g_2}Y)\,|\,\mathfrak{g}_i\in\mathfrak{G},\;i=1,2\}.
\end{equation*}
Depending on the choice of the group $ \mathfrak{G} $, we get  different classification of elements of the space $\conv\Bbb R^n$.
If, for example, $ \mathfrak{G} = GL(\Bbb R ^ n)\setminus \Bbb R ^ n $ is a general affine group of the space $ \Bbb R^ n $, then we obtain a more rough classification and if $ \mathfrak {G} = \Bbb R ^ n $ is a group of translations of the space $ \Bbb R ^ n $ then we obtain a thinner classification.

In this paper,   $ \mathfrak{G} $ is a group of homotheties namely  the semidirect product of the group of dilations and group of translations of the space $ \Bbb R^n $.

We introduce a subset
\begin{equation*}
\mathfrak C=\{X\in\conv\Bbb R^n\,|\,\,\intt X\ne\emptyset\}.
\end{equation*}
If $X\in\mathfrak{C}$, then we set $\wt X\overset{\texttt{def}}=\frac{X}{\sqrt[n]{V[X]}}$, where $V[X]$ is a volume of a convex compact $X$. Then $\Or_{\mathfrak{G}}(X)=\Or_{\mathfrak{G}}(\wt X)$ and if
$Y\in\mathfrak{C}$, then we get
\begin{equation*}
\rho[\Or_{\mathfrak{G}}(X),\Or_{\mathfrak{G}}(Y)]=\inf\{d_H(\wt X,\wt Y+\bold x)\,|\,\bold x\in\Bbb R^n\}=\rho[\Or_{\mathfrak{\Bbb R^n}}(X),\Or_{\mathfrak{\Bbb R^n}}(Y)].
\end{equation*}

In the space $\conv\Bbb R^n$ we consider the Cauchy problem for   SDEs
\begin{equation}\label{2.1}
D_HX(t)=\bold AX(t),\quad X(0)=X_0,\quad X_0\in\mathfrak C,
\end{equation}
where $D_H$ is a  Hukuhara derivative operator, $X(t)\in\conv\Bbb R^n$, $t\in\Bbb R_+$, $\bold A\in L(\Bbb R^n)$ is an orthogonal operator, that is $\bold A^*\bold A^{-1}=\bold I$, $\bold A^*$ is the adjoint operator.

Since $X_0\in\mathfrak{C}$ we have $X(t)\in\mathfrak{C}$ for all $t\ge 0$.

Let $X^*(t)$ be a program solution of the Cauchy problem \eqref{2.1} with the initial value $X^*(0)=X_0^*\in\mathfrak C$.

Next we give the definition of Lyapunov stability of solutions of the Cauchy problem \eqref{2.1}.

\begin{definition} Program solution $X^*(t)$ is said to be

1)    Lyapunov stable if for any $\e>0$ there exists a positive number $\dl=\dl(\e,X_0^*)$  such that, for all $X_0\in\mathfrak C$ the condition $\rho[\Or_{\mathfrak{G}}(X_0),\Or_{\mathfrak{G}}(X_0^*)]<\dl$ implies $\rho[\Or_{\mathfrak{G}}(X(t)),\Or_{\mathfrak{G}}(X^*(t))]<\e$ for all $t\ge 0$;

2) conditionally asymptotically stable with respect to the set $\mathfrak{M}\subset\mathfrak C$, if it is stable and  there exists a scalar
 $\s_0>0$ such that, for all $X_0\in\mathfrak{M}$ if $\rho[\Or_{\mathfrak{G}}(X_0),\Or_{\mathfrak{G}}(X_0^*)]<\s_0$  then we have
\begin{equation*}
\lim\limits_{t\to\infty}\rho[\Or_{\mathfrak{G}}(X(t)),\Or_{\mathfrak{G}}(X^*(t))]=0.
\end{equation*}
\end{definition}

In this paper, we investigate the stability and asymptotic stability of solutions of the Cauchy problem \eqref {2.1} in the sense of the above definition.

\section{Auxiliary results} For  $X,\,Y\in\mathfrak{C}$  we define the functional
\begin{equation*}
\Dl[X,Y]=\frac{V_1^n[X,Y]}{V^{n-1}[X]V[Y]}-1.
\end{equation*}
Here $V_1[X,Y]$ is the first mixed volume of convex compacts $X$ and $Y$.
Based on Brunn--Minkowski inequality \cite{bon}, we get $\Dl[X,Y]\ge 0$, and $\Dl[X,Y]=0$ if and only if $Y\in \Or_{\mathfrak{G}}(X)$. It is obvious that the functional $ \Dl[X, Y] $ depends only on the orbits $\Or_{\mathfrak{G}}(X)$ and $\Or_{\mathfrak{G}}(Y)$. This functional will be used as an analogue of the Lyapunov function.

The study of stability of solutions of the Cauchy problem \eqref {2.1} is based on the Theorem 1.2, proved in \cite{disk} on the basis of more accurate Brunn--Minkowski inequality established by  V.I. Discant.

For $X\in\mathfrak{C}$, let $R_{X}$,  $r_{X}$  be  the radii of the circumscribed and inscribed  balls for   convex compact  $ X $  respectively.  Let $X(t)$ and $Y(t)$ be one-parameter families of sets from $\mathfrak{C}$. It should be noted also that     $d_H((\widetilde{X})',\widetilde{Y})=\rho(\Or_{\mathfrak{G}}(X),\Or_{\mathfrak{G}}(Y))$, where $X'$ is a shift of set $X \in \conv \mathbb{R}^n$.
\begin{lemma}
Assume that
\begin{equation*}
R=\sup\limits_{t\ge 0}\{R_{\wt X(t)},R_{\wt Y(t)}\}<\infty,\quad r=\inf\limits_{t\ge 0}\{r_{\wt X(t)},r_{\wt Y(t)}\}>0.
\end{equation*}
Then there exist the positive constants $\e_0$, $C_1$ and $C_2$ which depend only on $ n $, $ R $ and $ r $ such that
\begin{equation}\label{2.25}
\Dl[X,Y]\le C_2\varrho[\Or_{\mathfrak{G}}(X),\Or_{\mathfrak{G}}(Y)]
\end{equation}
and if $\Dl[X,Y]<\e_0$  we have
\begin{equation}\label{2.25*}
C_1\varrho^{n^2}[\Or_{\mathfrak{G}}(X),\Or_{\mathfrak{G}}(Y)]\le \Dl[X,Y].
\end{equation}
\end{lemma}
{\it Proof.} Inequality \eqref{2.25*} is a direct consequence of the reasoning in the proof of Theorem 1.2 from  \cite{disk}. Next, we will prove  the inequality \eqref{2.25}. By definition of the metric, we get
\begin{equation*}
\gathered
(\wt Y(t))^{\prime}\subset\wt X(t)+\varrho\overline{B}_1(0),\quad
\varrho=\varrho[\Or_{\mathfrak{G}}(X(t)),\Or_{\mathfrak{G}}(Y(t))]=d_H(X(t),(Y(t))^{\prime}),
\endgathered
\end{equation*}
where $B_1(0) \subseteq \mathbb{R}^n$ is  an open unit ball with center at $x=0$.
By monotony of functional $V_1[X,Y]$ and using Brunn-Minkowski inequality,    we obtain
\begin{equation*}
\gathered
1\le V_1[\wt X(t),\wt Y(t)]=V_1[\wt X(t),(\wt Y(t))^{\prime}]\le V_1[\wt X(t),\wt X(t)+\varrho\overline{B}_1(0)]\\=
1+\varrho V_1[\wt X(t),\overline{B}_1(0)]
\le 1+\varrho V_1[R_{\wt X(t)}\overline{B}_1(0),\overline{B}_1(0)]\\=
1+\varrho R_{\wt X(t)}\upsilon\le 1+\varrho[\Or_{\mathfrak{G}}(\wt X(t)),\Or_{\mathfrak{G}}(\wt Y(t))]R\upsilon,
\endgathered
\end{equation*}
where $\upsilon=V[\overline{B}_1(0)]$. Thus, we get the equality
\begin{equation*}
\Dl[X(t),Y(t)]=V_1[\wt X(t),\wt Y(t)]-1\le C_2 \varrho[\Or_{\mathfrak{G}}(X(t)),\Or_{\mathfrak{G}}(Y(t))],\quad C_2=R\upsilon.
\end{equation*}
The lemma is proved.

In order to estimate  the changes of functionals $ V[X] $ and $ V_1[X, X^*] $ along the solutions of the Cauchy problem, we shall use the comparison method \cite{lak}.

Let $\bold k=\{k_1,\dots,k_{n-1}\}$  be a certain unordered set of indices, where $k_i\geq 0$,  and $ \Bbb K $ be the set of all such index sets.


Define the auxiliary functionals
\begin{equation*}
 \Xi_{\bold k}[X,X^*]=V[\bold A^{k_1}X,\bold A^{k_2}X,\dots,\bold A^{k_{n-1}}X,X^*]
\end{equation*}
and the functions $\zeta_{\bold k}(t)=\Xi_{\bold k}[X(t),X^*(t)]$.

Using the continuity of the functional of mixed volume, it is easy to show that for the sets  $\bold k=(k_1,\dots,k_p,0,\dots,0)$, $k_j\ge 1$, $j=1,\dots,p$ we have the formula
\begin{equation}\label{SS1}
\frac{d\zeta_{\bold k}(t)}{dt}=\zeta_{k_1+1,k_2,\dots,k_p,0,\dots,0}(t)+....+\zeta_{k_1,k_2,\dots,k_{p}+1,0,\dots,0}(t)+
(n-p)\zeta_{k_1,k_2,\dots,k_{p},1,0,\dots,0}(t),
\end{equation}
and for the sets $\bold k\in\Bbb K$ in which $k_j\ge 1$, $j=1$, $\dots$, $n-1$, we have
\begin{equation}\label{SS2}
\frac{d\zeta_{\bold k}(t)}{dt}=\zeta_{k_1+1,k_2,\dots,k_{n-1}}(t)+....+\zeta_{k_1,k_2,\dots,k_{n-1}+1}(t)+\zeta_{k_1-1,k_2-1,\dots,k_{n-1}-1}(t).
\end{equation}

Define the set
\begin{equation}
l_{\infty}=\{\{x_{\bold k}\}_{\bold k\in\Bbb K}\,\,|\,\,\sup\limits_{\bold k\in\Bbb K}|x_{\bold k}|<\infty\}
\end{equation}
and the norm $\|x\|_{l_{\infty}}=\sup\limits_{\bold k\in\Bbb K}|x_{\bold k}|$.

On the set $l_{\infty}$  the operations of addition  and  nonnegative scalar multiplication    are defined in a natural way.

 It is easy to see  that   $(l_{\infty},\|.\|_{\infty})$ is a Banach space.  Let us show that
 \{$\Xi_{\bold k}[X,X^*]\}_{\bold k\in\Bbb K}\in l_{\infty}$. In fact,
\begin{equation*}
\gathered
X^{\prime}\subset R_X\overline{B}_1(0),\quad (X^*)^{\prime}\subset R_{X^*}\overline{B}_1(0),\\
\Xi_{\bold k}[X,X^*]=\Xi_{\bold k}[X^{\prime},(X^*)^{\prime}]=V[\bold A^{k_1}X^{\prime},\bold A^{k_2}X^{\prime},\dots,\bold A^{k_{n-1}}X^{\prime},(X^*)^{\prime}]\\
\le R_{X}^{n-1}R_{X^*}V[\bold A^{k_1}\overline{B}_1(0),\bold A^{k_2}\overline{B}_1(0),\dots,\bold A^{k_{n-1}}\overline{B}_1(0),\overline{B}_1(0)]\le R_{X}^{n-1}R_{X^*}\upsilon.
\endgathered
\end{equation*}

Therefore, differential equations \eqref {SS1} and \eqref {SS2} can be represented in an abstract form
\begin{equation*}
\frac{d\zeta}{dt}=\Omega\zeta,
\end{equation*}
where $\zeta \in l_{\infty}$, $\Omega \colon l_{\infty}\to l_{\infty} $ is a linear operator. It is obvious that $\Omega \in L(l_{\infty})$ and $\|\Omega\|_{L(l_{\infty})}=n$.

Hence it follows that
\begin{equation}\label{MV}
V_1[X(t),X^*(t)]=\sum\limits_{\bold k\in\Bbb K}a_{\bold k}\Xi_{\bold k}[X_0,X_0^*],
\end{equation}
where $a_{\bold k}$ are the coefficients that do not depend on the operator $\bold A$.
 Operator $ \Omega $ is positive relative to the cone $l_{\infty}^+=\{\{x_{\bold k}\}_{\bold k\in\Bbb K}\in l_{\infty}\,\,|\,\,x_{\bold k}\ge 0\}$ and therefore the coefficients $a_{\bold k}$  are nonnegative.

Assuming that $X_0^*=X_0$ in formula \eqref{MV}, we get
\begin{equation}\label{MV1}
V[X(t)]=\sum\limits_{\bold k\in\Bbb K}a_{\bold k}M_{\bold k}[X_0],
\end{equation}
where $M_{\bold k}[X_0]=\Xi_{\bold k} [X_0,X_0]$.

Let us prove that
\begin{equation}\label{MV2}
\sum\limits_{\bold k\in\Bbb K}a_{\bold k}=e^{nt}.
\end{equation}
Let $X_0=\overline{B}_1(0)$, then as the    linear operator $A$ is orthogonal, we get $\bold A\overline{B}_1(0)=\overline{B}_1(0)$. It is obvious that $X(t)=e^t \overline{B}_1(0)$ is a unique solution of the Cauchy problem  (2.1) with initial conditions $X(0)=\overline{B}_1(0)$ and $M_\textbf{k}[\overline{B}_1(0)]=V[\overline B_1(0)]$. From (3.7) we obtain
$$
 e^{nt}V[\overline{B}_1(0)]=V[e^t\overline{B}_1(0)]=\sum\limits_{\bold k\in\Bbb K}a_{\textbf{k}}M_\textbf{k}[\overline{B}_1(0)]=\sum\limits_{\bold k\in\Bbb K}a_{\textbf{k}}V[\overline{B}_1(0)],
$$
and we get \eqref{MV2}.

\begin{lemma}
 For the volume $ V [X (t)] $  of solution $ X (t) $, $ X (0) = X_0 $  of linear differential equation \eqref {2.1}  the following estimate holds for all  $t \ge 0$
\begin{equation}
V[X_0]e^{nt}\le V[X(t)]\le M[X_0]e^{nt},
\end{equation}
where $M[X_0]=\max\limits_{\bold k\in\Bbb K}M_{\bold k}[X_0]$.
\end{lemma}

{\it Proof.} Applying the Brunn-Minkowski inequality, we obtain
\begin{equation*}
\frac{dV[X(t)]}{dt}=nV_1[X(t),\bold AX(t)]\ge nV[X(t)].
\end{equation*}
Hence  the inequality $V[X(t)]\ge e^{nt}V[X_0]$ is valid for $t\ge 0$.
\begin{equation*}
V[X(t)]\le \|\zeta(t)\|_{l_{\infty}}\le\|e^{t\Omega}\|_{L(l_{\infty})}\|\zeta_0\|_{l_{\infty}}\le
e^{t\|\Omega\|_{L(l_{\infty})}}M[X_0]=e^{nt}M[X_0].
\end{equation*}
The Lemma is proved.

\begin{lemma}
Assume  that $X(t)$, $X^*(t)$ are solutions  of differential equation \eqref{2.1}, $X(0)=X_0$, $\intt  X_0\ne\emptyset$, $X^*(0)=X_0^*$, $\intt  X^*_0\ne\emptyset$. Then there exist the positive constants $\e_0$,  $C_1$,  $C_2$   that depend on $X_0^*$, such that, from inequality
\begin{equation*}
\varrho[\Or_{\mathfrak{G}}(X_0),\Or_{\mathfrak{G}}(X_0^*)]<\s_0
\end{equation*}
 for all $t\ge 0$ we have
\begin{equation*}
\Dl[X(t),X^*(t)]\le C_2\varrho[\Or_{\mathfrak{G}}(X(t)),\Or_{\mathfrak{G}}(X^*(t))].
\end{equation*}
Moreover, for  $t\ge 0$ if $\Dl[X(t),X^*(t)]\le\e_0$, then we have
\begin{equation*}
C_1\varrho^{n^2}[\Or_{\mathfrak{G}}(X(t)),\Or_{\mathfrak{G}}(X^*(t))]\le\Dl[X(t),X^*(t)].
\end{equation*}
\end{lemma}
{\it Proof.} As a result of the assertion of Lemma 3.1, it suffices to prove that
\begin{equation*}
R=\sup\limits_{t\ge 0}\{R_{\wt X(t)},R_{\wt X^*(t)}\}\le R(X_0^*)<\infty,\quad
r=\inf\limits_{t\ge 0}\{r_{\wt X(t)},r_{\wt X^*(t)}\}\ge r(X_0^*)>0,
\end{equation*}
for all $X_0$ for which $\rho[\Or_{\mathfrak{G}}(X_0),\Or_{\mathfrak{G}}(X_0^*)]\le\s_0$, $\s_0$ is a sufficiently small positive constant.

Let $T$ be any positive number. For the Cauchy problem \eqref {2.1},  we consider the successive approximations for $t\in[0,T]$
\begin{equation*}
\gathered
X_0(t)=X_0,\quad X_m (t)=X_0+\int\limits_0^t\bold AX_{m-1}(s)\,ds.\\
\endgathered
\end{equation*}
Let $h_{X(t)}(p)$ be a support function for convex compact set $X(t)$.  Without loss of generality, we can assume that the origin of coordinates is at the center of inscribed ball, then  $h_{X_0}(p)\ge r_{X_0}$ for $p\in\partial B_1(0)$.  Next, we prove by  mathematical induction the inequality
\begin{equation}\label{N1}
\gathered
h_{X_m(t)}(p)\ge \sum\limits_{k=0}^m\frac{t^k}{k!}r_{X_0}.
\endgathered
\end{equation}
For $m=0$ this inequality is obvious. Suppose that it is true for $m=k-1$, then
\begin{equation}
\gathered
h_{X_k(t)}(p)=h_{X_0}(p)+\int\limits_0^th_{\bold AX_{k-1}(s)}(p)\,ds=
h_{X_0}(p)+\int\limits_0^th_{X_{k-1}(s)}(\bold A^*p)\,ds\\
\ge r_{X_0}+\int\limits_0^t
\sum\limits_{l=0}^{k-1}\frac{s^l}{l!}r_{X_0}\,ds=\sum\limits_{l=0}^{k}\frac{t^k}{k!}r_{X_0}.
\endgathered
\end{equation}
It is known \cite{lak}, that the successive approximations $ X_ {k} (t) $ converge uniformly with respect to $ t \in [0, T] $ to the solution $ X (t) $ of the Cauchy problem  \eqref{2.1} and therefore $\|h_{X_k(t)}-h_{X(t)}\|_{C(\partial B_1(0))}\to 0$ for $k\to\infty$. Thus, from inequality \eqref {N1}, it follows that $h_{X(t)}\ge e^tr_{X_0}$ for all $t\in\Bbb R_+$. Hence,
\begin{equation*}
\inf\limits_{t\ge 0}r_{\wt X(t)}=\inf\limits_{t\ge 0}\frac{r_{X(t)}}{\sqrt[n]{V[X(t)]}}\ge\inf\limits_{t\ge 0}\frac{e^{t}r_{X_0}}{e^t\sqrt[n]{M[X_0]}}=\frac{r_{X_0}}{\sqrt[n]{M[X_0]}}.
\end{equation*}
Similarly, we can show that $h_{X(t)}\le e^tR_{X_0}$ for all $t\ge 0$, and therefore
\begin{equation*}
\sup\limits_{t\ge 0}R_{\wt X(t)}=\sup\limits_{t\ge 0}\frac{R_{X(t)}}{\sqrt[n]{V[X(t)]}}\le\sup\limits_{t\ge 0}\frac{e^{t}R_{X_0}}{e^t\sqrt[n]{V[X_0]}}=\frac{R_{X_0}}{\sqrt[n]{V[X_0]}}.
\end{equation*}
Thus,
\begin{equation*}
r\ge\min\Big[\frac{r_{X_0}}{\sqrt[n]{M[X_0]}}
,\frac{r_{X_0^*}}{\sqrt[n]{M[X_0^*]}}\Big]
\end{equation*}
and
\begin{equation*}
R\le\max\Big[\frac{R_{X_0}}{\sqrt[n]{V[X_0]}}
,\frac{R_{X_0^*}}{\sqrt[n]{V[X_0^*]}}\Big].
\end{equation*}
By the continuity of the functionals $ R_X $, $ r_X $, $ V [X] $ and $ M [X] $, there exists a positive constant $\e_0<\s_0$ such that, the inequality $\rho[\Or_{\mathfrak{G}}(X_0),\Or_{\mathfrak{G}}(X_0^*)]<\e_0$ implies the estimates
\begin{equation*}
\gathered
\Big|\frac{R_{X_0}}{\sqrt[n]{V[X_0]}}-\frac{R_{X_0^*}}{\sqrt[n]{V[X_0^*]}}\Big|<\frac{R_{X_0^*}}{2\sqrt[n]{V[X_0^*]}},\\
\Big|\frac{r_{X_0}}{\sqrt[n]{M[X_0]}}-\frac{r_{X_0^*}}{\sqrt[n]{M[X_0^*]}}\Big|<\frac{r_{X_0^*}}{2\sqrt[n]{M[X_0^*]}}.
\endgathered
\end{equation*}
So, we get
\begin{equation*}
R\le\frac{3R_{X_0^*}}{2\sqrt[n]{V[X_0^*]}},\quad r\ge\frac{r_{X_0^*}}{2\sqrt[n]{M[X_0^*]}}.
\end{equation*}
The Lemma is proved.

Consider the particular case of the Cauchy problem \eqref {2.1} when $ n = 2 $ and for some positive integer $ m $ the equality $\bold A^m=\bold I$ is valid.

\begin{lemma} Assume that $X(t)$ and $X^*(t)$ are solutions of the Cauchy problem  \eqref{2.1}  with initial conditions $X(0)=X_0$, $X^*(0)=X_0^*$.  Then for odd  $m$, $m\geq3$  we have the formula
\begin{equation*}
\gathered
S[X(t),X^*(t)]=\frac{1}{m}(e^{2t}+2\sum\limits_{q=1}^{[m/2]}e^{2t\cos\frac{2\pi q}{m}})S[X_0,X_0^*]\\+\frac{1}{m^2}\sum\limits_{p=1}^{m-1}\Big((m-p)e^{2t}+2\sum\limits_{q=1}^{[m/2]}[(m-p)\cos\frac{2\pi pq}{m}\\+2t\sin\frac{2\pi pq}{m}\sin\frac{2\pi q}{m}]e^{2t\cos\frac{2\pi q}{m}}\Big)(S[X_0,\bold A^p X_0^*]+S[X_0^*,\bold A^p X_0]).
\endgathered
\end{equation*}
For even  $m$, $m\geq4$ we have the formula
\begin{equation*}
\gathered
S[X(t),X^*(t)]=\frac{1}{m}(e^{2t}+e^{-2t}+2\sum\limits_{q=1}^{(m-2)/2}e^{2t\cos\frac{2\pi q}{m}})S[X_0]\\+\frac{1}{m^2}\sum\limits_{p=1}^{m-1}\Big((m-p)(e^{2t}+(-1)^pe^{-2t})+2\sum\limits_{q=1}^{(m-2)/2}[(m-p)\cos\frac{2\pi pq}{m}\\+2t\sin\frac{2\pi pq}{m}\sin\frac{2\pi q}{m}]e^{2t\cos\frac{2\pi q}{m}}\Big)(S[X_0,\bold A^p X_0^*]+S[X_0^*,\bold A^p X_0]).
\endgathered
\end{equation*}
\end{lemma}
{\it Proof.} Let $S[X,Y]$ be a functional  of  Minkowski mixed  area,  then  the auxiliary functions
\begin{equation*}
\gathered
\xi_k(t)=\frac{1}{2}(S[X(t),\bold A^kX^*(t)]+S[X^*(t),\bold A^kX(t)]), \quad k=0, \dots, m-1,
\endgathered
\end{equation*}
 satisfy the system of differential equations of $ m$ -th order (comparison system)
\begin{equation*}
\frac{d\xi}{dt}=\Omega\xi(t),
\end{equation*}
where $\xi(t)\in\Bbb R^m$, $\Omega\in\Bbb R^{m\times m}$ is the matrix, the non-zero elements of which have the form $\omega_{12}=2$, $\omega_{ij}=\omega_{m1}=1$, $|i-j|=1$, $(i,j)\ne(1,2)$.

It is known  \cite{krein} that the solution of comparison system  has the form
\begin{equation}\label{RF}
\gathered
\xi(t)=-\frac{1}{2\pi i}\oint\limits_{\Gamma}e^{\lm t}R_{\Omega}(\lm)\,d\lm\xi(0).
\endgathered
\end{equation}
Here $R_{\Omega}(\lm)$ is the
 resolvent of matrix $\Omega$, $\Gamma$ is the circuit consisting of a finite number of closed Jordan curves, oriented in the positive direction, covering the  spectrum of $ \s (\Omega) $ of  matrix $\Omega$.

Next, we find the spectrum  $\s(\Omega)$ and resolvent $R_{\Omega}(\lm)$ of matrix $\Omega$.
Let $f=(f_0,\dots,f_{m-1})\in\Bbb C^m$, $x=(x_0,\dots,x_{m-1})\in\Bbb C^m$ and consider the linear equation
\begin{equation}\label{HU}
(\Omega-\lm E)x=f.
\end{equation}
Equation \eqref{HU} is equivalent to the boundary value problem for finite-difference equation of second order
\begin{equation*}
x_{k-1}+x_{k+1}-\lm x_k=f_k,\quad k=1,\dots,m-1
\end{equation*}
with the boundary conditions
\begin{equation*}
2x_1-\lm x_0=f_0,\quad x_m=x_0.
\end{equation*}
If $f=0$ and equation \eqref{HU} has only the trivial solution, then $\lm\in\varrho(\Omega)$, where
$\varrho(\Omega)$ is the resolvent set of  matrix $\Omega$. Moreover, $\s(\Omega)=\Bbb C\setminus\varrho(\Omega)$. The general solution of the homogeneous difference equation has the form
\begin{equation*}
x_k=c_1q_1^k+c_2q_2^k,
\end{equation*}
where $q_1$ and $q_2$  are the roots of a quadratic equation $q^2+1=\lm q$, $c_1$ and $c_2$  are an arbitrary constants.  From the boundary conditions it follows that
\begin{equation*}
\gathered
2x_1-\lm x_0=(c_1-c_2)(q_1-q_2)=0,\\
c_1q_1^m+c_2q_2^m=x_0.
\endgathered
\end{equation*}

There are two cases: $c_1=c_2$ or $q_1=q_2$. In the first case $c_1(q_1^m+q_2^m-2)=0$, and if $(q_1^m+q_2^m-2)\ne 0$, then there is only the trivial solution of the equation \eqref{HU}. In the second case, the condition $ q_1 = q_2 $ implies that $q_1=q_2=1$ or $q_1=q_2=-1$. If $q_1=q_2=1$, then $x_k=c$  is a solution of \eqref{HU} for any $c$, i.e., $2\notin\varrho(\Omega)$. If
$q_1=q_2=-1$, then $x_k=(-1)^k(c_1+c_2)$, and for  $k=m$ we get $(1+(-1)^{m+1})(c_1+c_2)=0$, and if $(1+(-1)^{m+1})\ne 0$, then there is only the trivial solution of the equation \eqref{HU}. Thus, if  $m$ is the odd number then $-2\in\varrho(\Omega)$, otherwise it is obvious that $-2\notin\varrho(\Omega)$.

So,  for matrix spectrum we  have
\begin{equation*}
\s(\Omega)=\Big\{2\cos\frac{2\pi q}{m}\,|\quad q=0,1,...,\Big[\frac{m}{2}\Big]\Big\}.
\end{equation*}

 Next, we consider the equation \eqref{HU} in the general case when $f\ne 0$. It is easy to show that the general solution of the inhomogeneous finite-difference equation has the form
\begin{equation*}
x_k=\begin{cases}
c_1+c_2,\quad k=0,\\
c_1q_1^k+c_2q_2^k+\sum\limits_{p=0}^{k-1}\frac{q_1^{k-p}-q_2^{k-p}}{q_1-q_2}f_p,\quad k\ge 1.
\end{cases}
\end{equation*}
Taking into account the boundary conditions,
\begin{equation*}
2x_1-\lm x_0=f_0,\quad x_m=x_0
\end{equation*}
we obtain:
\begin{equation*}
\gathered
c_1-c_2=\frac{f_0}{q_2-q_1},\\
(1-q_1^m)c_1+(1-q_2^m)c_2=\frac{q_1^m-q_2^m}{q_1-q_2}f_0+\sum\limits_{p=1}^{m-1}\frac{q_1^{m-p}-q_2^{m-p}}{q_1-q_2}f_p.
\endgathered
\end{equation*}
Since $x_0=c_1+c_2$, we get
\begin{equation*}
x_0=\frac{q_1^m-q_2^m}{(q_1-q_2)(2-q_1^m-q_2^m)}f_0+2\sum\limits_{p=1}^{m-1}\frac{q_1^{m-p}-q_2^{m-p}}{(q_1-q_2)(2-q_1^m-q_2^m)}f_p.
\end{equation*}
Similarly,  we can find  all $ x_k $, $ k = 1, \dots, m-1 $.

 For $q_i=q_i(\lm)$, $i=1,2$ from \eqref{RF}, we can get the equality for mixed area
$S[X(t),X^*(t)]$ of solution  of the Cauchy problem  \eqref{2.1}
\begin{equation}\label{SQ}
\gathered
S[X(t),X^*(t)]=-\frac{1}{2\pi i}\oint\limits_{\Gamma}\frac{(q_1^m(\lm)-q_2^m(\lm))e^{\lm t}}{(q_1(\lm)-q_2(\lm))(2-q_1^m(\lm)-q_2^m(\lm))}\,d\lm\,\,S[X_0,X_0^*]\\-
\sum\limits_{p=1}^{m-1}\frac{1}{2\pi i}\oint\limits_{\Gamma}\frac{(q_1^{m-p}(\lm)-q_2^{m-p}(\lm))e^{\lm t}}{(q_1(\lm)-q_2(\lm))(2-q_1^m(\lm)-q_2^m(\lm))}\,d\lm\,\,(S[X_0^*,\bold A^pX_0]+S[X_0,\bold A^pX_0^*]).
\endgathered
\end{equation}
Thus, further calculations are reduced to finding the corresponding integrals in formula \eqref{SQ}.  The contour of integration can be represented as
\begin{equation*}
\Gamma=\bigcup\limits_{q=0}^{[m/2]}\Gamma_q,\quad \Gamma_q=\{\lm\in\Bbb C\,:\,\,\Big|\lm-2\cos\frac{2\pi q}{m}\Big|=\e\}, q=0,\dots,\Big[\frac{m}{2}\Big],
\end{equation*}
where $\e$ is a sufficiently small positive number.

Consider the integral $I_q=\oint\limits_{\Gamma_q}\frac{(q_1^m(\lm)-q_2^m(\lm))e^{\lm t}}{(q_1(\lm)-q_2(\lm))(2-q_1^m(\lm)-q_2^m(\lm))}\,d\lm$ for $q\ne 0,m/2$.
Then $\lm=2\cos\frac{2\pi q}{m}+\e e^{i\varphi}$, $\varphi\in[0,2\pi]$.
In this case, for sufficiently small $ \e $ the small parameter power series expansions   are valid:
\begin{equation*}
\gathered
2-(q_1^m(\lm)+q_2^m(\lm))=-\frac{m^2\e^2(\sin\varphi-i\cos\varphi)^2}{4\sin^2\frac{2\pi q}{m}}+o(\e^2),\\
q_1(\lm)-q_2(\lm)=2i\sin\frac{2\pi q}{m}+o(1).
\endgathered
\end{equation*}
So, we get
\begin{equation}\label{I0}
I_q=\int\limits_0^{2\pi}\frac{2\frac{m\e(\sin\varphi-i\cos\varphi)}{2\sin\frac{2\pi q}{m}}\e i e^{i\varphi}e^{2t\cos\frac{2\pi q}{m}}}{-\frac{m^2\e^2(\sin\varphi-i\cos\varphi)^2}{4\sin^2\frac{2\pi q}{m}}2i\sin\frac{2\pi q}{m}}\,d\varphi+o(1)=-\frac{4\pi i}{m}e^{2t\cos\frac{2\pi q}{m}}+o(1).
\end{equation}
Similarly, we obtain
\begin{equation}\label{I1}
\gathered
J_q=\oint\limits_{\Gamma_q}\frac{(q_1^{m-p}(\lm)-q_2^{m-p}(\lm))e^{\lm t}}{(q_1(\lm)-q_2(\lm))(2-q_1^m(\lm)-q_2^m(\lm))}\,d\lm=
-\frac{4i\pi(m-p)\cos\frac{2\pi pq}{m}}{m^2}e^{2t\cos\frac{2\pi q}{m}}\\-\frac{8\pi i}{m^2}t\sin\frac{2\pi pq}{m}\sin\frac{2\pi q}{m}e^{2t\cos\frac{2\pi q}{m}}+o(1),\quad
q\ne 0,\quad q\ne m/2.
\endgathered
\end{equation}

If $q=0$, then
\begin{equation}\label{I2}
I_0=e^{2t}\int\limits_0^{2\pi}\frac{2me^{i\varphi/2}\sqrt{\e} i\e e^{i\varphi}}{2e^{i\varphi/2}\sqrt{\e}(-m^2\e e^{i\varphi})}\,d\varphi+o(1)=-\frac{2\pi i e^{2t}}{m}+o(1).
\end{equation}
Similarly, we obtain
\begin{equation}\label{I3}
J_0=e^{2t}\int\limits_{0}^{2\pi}\frac{2(m-p)e^{i\varphi/2}\sqrt{\e}i\e e^{i\varphi}}{2e^{i\varphi}\sqrt{\e}(-m^2\e e^{i\varphi})}\,d\varphi+o(1)=-\frac{2\pi i(m-p)}{m^2}e^{2t}+o(1).
\end{equation}
If the number $ m $ is even, then  it is necessary to calculate the integrals $I_{m/2}$ and $J_{m/2}$. In this case we get
\begin{equation}\label{I4}
I_{m/2}=-\frac{2\pi i}{m}e^{-2t}+o(1).
\end{equation}
\begin{equation}\label{I5}
J_{m/2}=(-1)^p\frac{2\pi i(m-p)}{m^2}e^{-2t}+o(1).
\end{equation}

Substituting the calculated integrals \eqref{I0}---\eqref{I5} in \eqref{RF},  we obtain the assertion of Lemma. This completes the proof.

\section{Main result}  In this section, we aim to establish stability conditions of  the program solutions of  differential equation \eqref{2.1}.

\begin{theorem} Assume that   $\bold A$  is an orthogonal operator, then any program solution  $X^*(t)$, $X^*(0)=X_0^*$, $\intt X_0^*\ne\emptyset$ is Lyapunov stable.
\end{theorem}
{\it Proof.} Define the function
\begin{equation*}
\varphi(t)=\Dl[X(t),X^*(t)]=\frac{V_1^{n}[X(t),X^*(t)]-V^{n-1}[X(t)]V[X^*(t)]}{V^{n-1}[X(t)]V[X^*(t)]}.
\end{equation*}
 By formulas \eqref{MV} and \eqref{MV1},  the numerator of $\varphi(t)$ can be represented as
\begin{equation*}
\gathered
V_1^{n}[X(t),X^*(t)]-V^{n-1}[X(t)]V[X^*(t)]\\=\Big(\sum\limits_{\bold k\in\Bbb K}a_{\bold k}\Xi_{\bold k}[X_0,X_0^*]\Big)^n-\Big(\sum\limits_{\bold k\in\Bbb K}a_{\bold k}M_{\bold k}[X_0]\Big)^{n-1}
\sum\limits_{\bold k\in\Bbb K}a_{\bold k}M_{\bold k}[X_0^*]\\=\sum\limits_{\bold k_1\in\Bbb K}...
\sum\limits_{\bold k_n\in\Bbb K}(\Xi_{\bold k_1}[X_0,X_0^*]...\Xi_{\bold k_n}[X_0,X_0^*]-M_{\bold k_1}[X_0]\dots M_{\bold k_{n-1}}[X_0]M_{\bold k_n}[X_0^*])a_{\bold k_1}\dots a_{\bold k_n}.
\endgathered
\end{equation*}
Therefore, taking into account Lemma 3.2, we obtain the estimate
\begin{equation*}
\gathered
\varphi(t)\le e^{-n^2t}\sum\limits_{\bold k_1\in\Bbb K}...
\sum\limits_{\bold k_n\in\Bbb K}(\Xi_{\bold k_1}[\wt X_0,\wt X_0^*]...\Xi_{\bold k_n}[\wt X_0,\wt X_0^*]\\-M_{\bold k_1}[\wt X_0]\dots M_{\bold k_{n-1}}[\wt X_0]M_{\bold k_n}[\wt X_0^*])a_{\bold k_1}\dots a_{\bold k_n}.
\endgathered
\end{equation*}
By the definition of the metric,
\begin{equation*}
\gathered
(\wt X_0)^{\prime}\subset\wt X_0^*+\rho_0B_1(0),\quad \rho_0=\rho[\Or_{\mathfrak{G}}(\wt X_0),\Or_{\mathfrak{G}}(\wt X_0^*)].
\endgathered
\end{equation*}
By the monotony of the mixed volume functional, we obtain
\begin{equation*}
\gathered
\Xi_{\bold k}[\wt X_0,\wt X_0^*]=\Xi_{\bold k}[(\wt X_0)^{\prime},\wt X_0^*]\le
\Xi_{\bold k}[\wt X_0^*,\wt X_0^*]+\sum\limits_{k=1}^{n-1}C_{n-1}^k\rho_0^kR^{n-k}_{\wt X_0^*}\upsilon\\
=M_{\bold k}[\wt X_0^*]+\sum\limits_{k=1}^{n-1}C_{n-1}^k\rho^k_0R^{n-k}_{\wt X_0^*}\upsilon.
\endgathered
\end{equation*}
From the inclusion
\begin{equation*}
\gathered
(\wt X_0^*)^{\prime}\subset\wt X_0+\rho_0\overline{B}_1(0)
\endgathered
\end{equation*}
and the monotony of the mixed volume, it follows the inequality
\begin{equation*}
\gathered
M_{\bold k}[\wt X_0^*]\le M_{\bold k}[\wt X_0]+\sum\limits_{k=1}^{n-1}C_{n-1}^k\rho^k_0(R_{\wt X_0^*}+\rho_0)^{n-k}\upsilon.
\endgathered
\end{equation*}
Choose a positive number $ \e_1 $ so that for all $ \rho $, $ 0 <\rho <\e_1 $  the following inequalities hold
\begin{equation*}
\gathered
\sum\limits_{k=1}^{n-1}C_{n-1}^k\rho^k(R_{\wt X_0^*}+\rho)^{n-k}\upsilon\le 2(n-1)\rho R^{n-1}_{\wt X_0^*}\upsilon,\\
\sum\limits_{k=1}^{n-1}C_{n-1}^k\rho^kR^{n-k}_{\wt X_0^*}\upsilon\le 2(n-1)\rho R^{n-1}_{\wt X_0^*}\upsilon,\\
\e_1<\frac{M[X_0^*]}{2(n-1)R^{n-1}_{\wt X_0^*}\upsilon}.
\endgathered
\end{equation*}
Then for all  $X_0$ such that $\rho_0<\e_1$ the following inequality holds
\begin{equation*}
\gathered
\varphi(t)\le e^{-n^2t}\sum\limits_{\bold k_1\in\Bbb K}...
\sum\limits_{\bold k_n\in\Bbb K}\Big[(M_{\bold k}[\wt X_0^*]+2(n-1)\rho_0 R^{n-1}_{\wt X_0^*}\upsilon)^n\\-(M_{\bold k}[\wt X_0^*]-2(n-1)\rho_0 R^{n-1}_{\wt X_0^*}\upsilon)^{n-1}M_{\bold k}[\wt X_0^*]\Big]a_{\bold k_1}\dots a_{\bold k_n}.
\endgathered
\end{equation*}
Applying  Lagrange's theorem on finite increments  for function
\begin{equation*}
\gathered
f(\rho)=(M_{\bold k}[\wt X_0^*]+2(n-1)\rho R^{n-1}_{\wt X_0^*}\upsilon)^n\\-(M_{\bold k}[\wt X_0^*]-2(n-1)\rho R^{n-1}_{\wt X_0^*}\upsilon)^{n-1}M_{\bold k}[\wt X_0^*]
\endgathered
\end{equation*}
we get the following  estimate
\begin{equation*}
\gathered
\varphi(t)\le e^{-n^2t}A\rho_0\sum\limits_{\bold k_1\in\Bbb K}...
\sum\limits_{\bold k_n\in\Bbb K}a_{\bold k_1}\dots a_{\bold k_n}=
e^{-n^2t}A\rho_0\Big(\sum\limits_{\bold k\in\Bbb K}a_{\bold k}\Big)^n,
\endgathered
\end{equation*}
where
\begin{equation*}
\gathered
A=2(n-1)R_{\wt X_0^*}\upsilon(M[X_0^*]+2(n-1)R_{\wt X_0^*}^{n-1}\upsilon\e_1)^{n-1}
(M[X_0^*]\\+2n(n-1)R_{\wt X_0^*}^{n-1}\upsilon\e_1).
\endgathered
\end{equation*}
By formula \eqref{MV2}, we get
\begin{equation*}
\gathered
\varphi(t)\le A\rho_0.
\endgathered
\end{equation*}
From the assertion of Lemma 3.1 it follows that there exists $ \e_0> 0 $ such that for $\rho[ \Or_{\mathfrak{G}}(X_0), \Or_{\mathfrak{G}}(X_0^*)]<\e_m$,
$\e_m=\min[\e_0,\e_1]$ the estimate holds
\begin{equation*}
\gathered
\rho[\Or_{\mathfrak{G}}(X(t)),\Or_{\mathfrak{G}}(X^*(t))]\le\Big(\frac{C_2A}{C_1}\Big)^{1/n^2}\rho^{1/n^2}[\Or_{\mathfrak{G}}(X_0),\Or_{\mathfrak{G}}(X_0^*)].
\endgathered
\end{equation*}
For $\e>0$ we choose $\dl(\e)=\min\Big\{\e_m,\e^{n^2}\Big(\frac{C_2A}{C_1}\Big)^{-1/n^2}\Big\}$,  then for all  $t\ge 0$ the inequality
$\rho[\Or_{\mathfrak{G}}(X(t)),\Or_{\mathfrak{G}}(X^*(t))]<\e$ is fulfilled. This completes the proof.

\textbf{ Remark.} The assertion of Theorem 4.1 is valid  if  the orthogonality condition  for operator $ \bold A $  is replaced  by the condition $\sup\limits_{k\in\Bbb Z}\|\bold A^k\|<\infty$.

Indeed, in this case,  by  Theorem 6.1  about stable operators \cite{krein}, there is an orthogonal operator $ \bold A_1 $ and a nonsingular operator
$\bold T$ such that, $\bold A_1=\bold T^{-1}\bold A\bold T$. In the Cauchy problem \eqref {2.1} we make the change of variables $ X = \bold TY $, then this problem is of
the form
\begin{equation}\label{2.5bis}
D_HY(t)=\bold AY(t),\quad Y(0)=Y_0.
\end{equation}
Thus, it is obvious that $\Or_{\mathfrak{G}}(\bold TX)=\bold T\Or_{\mathfrak{G}}(X)$ and the stability  problem of  solution $ X^ * (t) $ of the Cauchy problem \eqref {2.1} is equivalent to the stability  problem of solution  $Y^*(t)=\bold T^{-1}X^*(t)$ of the Cauchy problem \eqref{2.5bis}.

%

\begin{theorem} Assume that $n=2$ and there exists a positive integer number $m$ such that, operator  $\bold A^m=\bold I$, then any  solution $ X^* (t) $ is conditional  Lyapunov asymptotically stable relative to the set
$$
\mathfrak{M}=\Big\{X_0\,\,|\,\intt X_0\ne\emptyset,\,\sum\limits_{p=0}^{m-1}\bold A^pX_0\in\Or_{\mathfrak{G}}\Big( \sum\limits_{p=0}^{m-1}\bold A^pX_0^*\Big)\Big\}.
$$
\end{theorem}
{\it Proof.} Stability  of solution $ X^* (t) $ is the consequence of Theorem 4.1. Let us prove the condition of attraction of $ X^* (t) $ relative to the set  $\mathfrak{M}$. Consider the function
$$
\varphi(t)=\Dl[X(t),X^*(t)]=\frac{S^2[X(t),X^*(t)]}{S[X(t)]S[X^*(t)]}-1.
$$
As a result of the assertion of Lemma 3.4 we get
$$
\gathered
S[X(t),X^*(t)]=e^{2t}\Big(\frac{1}{m}S[X_0,X_0^*]+\frac{1}{m^2}\sum\limits_{p=1}^{m-1}(m-p)(S[\bold A^pX_0,X_0^*]\\+S[X_0,\bold A^p X_0^*]\Big)+o(e^{2t})=\frac{e^{2t}}{m^2}S\Big[\sum\limits_{k=0}^{m-1}\bold A^kX_0,\sum\limits_{k=0}^{m-1}\bold A^kX_0^*\Big]+o(e^{2t}),\quad t\to\infty.
\endgathered
$$
$$
\gathered
S[X(t)]=e^{2t}\Big(\frac{1}{m}S[X_0]+\frac{2}{m^2}\sum\limits_{p=1}^{m-1}(m-p)S[\bold A^pX_0,X_0]\Big)+o(e^{2t})\\=\frac{e^{2t}}{m^2}S\Big[\sum\limits_{k=0}^{m-1}\bold A^kX_0\Big]+o(e^{2t}),\quad t\to\infty.
\endgathered
$$

Then we have
$$
\varphi(t)=\frac{\frac{e^{4t}}{m^4}\Big(S^2\Big[\sum\limits_{k=0}^{m-1}\bold A^kX_0,\sum\limits_{k=0}^{m-1}\bold A^kX_0^*\Big]-S\Big[\sum\limits_{k=0}^{m-1}\bold A^kX_0^*\Big]S\Big[\sum\limits_{k=0}^{m-1}\bold A^kX_0\Big]\Big)+o(e^{4t})}{\frac{e^{4t}}{m^4}S\Big[\sum\limits_{k=0}^{m-1}\bold A^kX_0^*\Big]S\Big[\sum\limits_{k=0}^{m-1}\bold A^kX_0\Big]+o(e^{4t})}.
$$
By Lemma 3.1 $\rho[\Or_{\mathfrak{G}}(X(t)),\Or_{\mathfrak{G}}(X^*(t))]\to 0$ for $t\to\infty$ if and only if when $\varphi(t)\to 0$ for $t\to\infty$. It's obvious that $\varphi(t)\to 0$ for $t\to\infty$ if and only if when
$$
S^2\Big[\sum\limits_{k=0}^{m-1}\bold A^kX_0,\sum\limits_{k=0}^{m-1}\bold A^kX_0^*\Big]-S\Big[\sum\limits_{k=0}^{m-1}\bold A^kX_0^*\Big]S\Big[\sum\limits_{k=0}^{m-1}\bold A^kX_0\Big]=0.
$$
By  Brunn--Minkowski theorem,  the last equality is valid if and only if when
\begin{equation*}
\sum\limits_{k=0}^{m-1}\bold A^kX_0\in \Or_{\mathfrak{G}}\Big(\sum\limits_{k=0}^{m-1}\bold A^kX_0^*\Big).
\end{equation*}
The Theorem is proved.

\section{Example} Assume that the operator $ \bold A $ is the rotation operator in the positive direction by the angle  $\frac{2\pi}{m}$. Then,  the matrix
$ A $ of the linear operator $ \bold A $ in the canonical basis has the form
\begin{equation*}
A=\begin{pmatrix}
\cos\frac{2\pi}{m}&-\sin\frac{2\pi}{m}\\
\sin\frac{2\pi}{m}&\cos\frac{2\pi}{m}
\end{pmatrix}.
\end{equation*}
Let $h_X(\bold p)$ be a  support function of a convex compact $X\in\conv\Bbb R^2$, $H_X(\theta)=h_X(\cos\theta,\sin\theta)$.
Then
\begin{equation*}
H_{\sum\limits_{k=0}^{m-1}\bold A^kX_0}(\theta)=\sum\limits_{k=0}^{m-1}H_{\bold A^kX_0}(\theta)=
\sum\limits_{k=0}^{m-1}H_{X_0}\Big(\theta-\frac{2\pi k}{m}\Big).
\end{equation*}
For function $H_{X_0}(\theta)$  we can obtain  the Fourier series expansion
\begin{equation*}
H_{X_0}(\theta)=\sum\limits_{p=-\infty}^{\infty}H_pe^{ip\theta}.
\end{equation*}
Hence, we obtain
\begin{equation*}
\gathered
H_{\sum\limits_{k=0}^{m-1}\bold A^kX_0}(\theta)=\sum\limits_{k=0}^{m-1}\sum\limits_{p=-\infty}^{\infty}H_pe^{ip(\theta-\frac{2\pi k}{m})}=
\sum\limits_{p=-\infty}^{\infty}H_pe^{ip\theta}\sum\limits_{k=0}^{m-1}e^{\frac{2\pi ipk}{m}}\\=
\sum\limits_{p=-\infty}^{\infty}H_{pm}e^{ipm\theta}.
\endgathered
\end{equation*}
By Theorem 4.2 from conditions
\begin{equation*}
\int\limits_0^{2\pi}(H_{X_0}(\theta)-H_{X_0^*}(\theta))e^{-mp\theta}d\theta=0, p\in\Bbb Z_+
\end{equation*}
it follows that
\begin{equation*}
\lim\limits_{t\to\infty}\rho[\Or_{\mathfrak{G}}(X(t)),\Or_{\mathfrak{G}}(X^*(t))]= 0,
\end{equation*}
provided that $\rho[\Or_{\mathfrak{G}}(X_0),\Or_{\mathfrak{G}}(X^*_0)]$ is the sufficiently small positive number.

Thus, for each solution of the Cauchy problem $ X^ * (t) $, $ X^ * (0) \in \mathfrak {C} $ there is an infinite dimensional variety of solutions
$ X (t) $, that are attracted to the solution  $X^*(t)$.

\section{Conclusion}  By Theorem 4.1 we can conclude that the solution of the Cauchy problem \eqref {2.1} has a stable form. Theorem 4.2 strengthens this result for the case of two-dimensional space and a periodic operator.  It suggests that for each solution there is an infinite-dimensional variety of solutions that are attracted to the program solution.  For further study  it is of interest to generalize the Theorem 4.2 for spaces of dimension greater than 2, and also for stable nonperiodic operators $ \bold A $. The main hypothesis concerning this case is that the forms of all solutions of the Cauchy problem \eqref {2.1}
asymptotically tend to  a ball shape.

\bigskip

\bigskip

This work was supported by grant of Ministry of Education and Science of Ukraine [grant number  0116U004691].

\bigskip

\bigskip

The author is grateful to the Prof. G.T. Bhaskar for discussions and valuable comments that improved the
manuscript.

\enddocument